\documentclass[a4paper, 12pt]{article}

\usepackage{setspace}
\usepackage[latin1]{inputenc}
\usepackage{graphics}
\usepackage{epsfig}
\usepackage{bbm}
\usepackage{wrapfig}
\usepackage{float}
\usepackage{setspace}
\usepackage{picins}
\usepackage{amssymb}
\usepackage{mathrsfs}
\usepackage{amsthm}
\usepackage{amsmath}
\usepackage{paralist}
\usepackage[T1]{fontenc}

\begin{document}

\normalsize

\newtheorem{theorem}{Theorem}
\newenvironment{void}{}

\newtheorem{proposition}[theorem]{Proposition}
\newtheorem{lemma}[theorem]{Lemma}
\newtheorem{corollary}[theorem]{Corollary}
\newtheorem{remark}[theorem]{Remark}
\newtheorem{example}[theorem]{Example}
\newtheorem{conjecture}[theorem]{Conjecture}
\newtheorem*{results}{Results}
\newtheorem{program}{Program}
\newtheorem*{notation}{Notation}
\newtheorem{algorithm}[theorem]{Algorithm}

\theoremstyle{definition}
\theoremstyle{remark}

\author{Mykhaylo Tyomkyn}
\title{An improved bound for the Manickam-Mikl\'os-Singhi conjecture}
\maketitle

\begin{abstract}
We show that for $n>k(4e\log k)^k$ every set $\{x_1,\cdots, x_n\}$ of $n$ real numbers with $\sum_{i=0}^{n}x_i \geq 0$ has at least $\binom{n-1}{k-1}$ $k$-element subsets of a non-negative sum. This is a substantial improvement on the best previously known bound of $n>(k-1)(k^k+k^2)+k$, proved by Manickam and Mikl\'os~\cite{MM} in 1987.
\end{abstract}

\section{Introduction}

Throughout this paper we shall use the standard abbreviation $[n]$ for the set of natural numbers from $1$ to $n$ and $\mathscr{P}([n])$ for its power set. For a set $C$ and an integer $k$ we write, as usual, $C^{(k)}$ for $\left\{A \subset C\colon \left|A\right|= k \right\}$. For a set $\mathcal{X}$ of real numbers we shall write $\sum \mathcal{X}$ for the sum of its elements.

Let $\mathcal{X} = \left\{x_1, x_2, \dotsc x_n\right\}$ be a collection of $n$ real numbers whose sum is non-negative. How many subsets of $\mathcal{X}$ can be guaranteed to have a non-negative sum? The answer is $2^{n-1}$, since for every $C\subset[n]$ at least one of $\sum_{i\in C} x_i$ and $\sum_{i\in [n] \setminus C} x_i$ is non-negative. The bound is tight, as can be easily seen by taking any collection of numbers whose total sum is $0$ but any partial sum is either strictly positive or negative, for instance $x_1 = n-1 $ and $x_2 = x_3 = \dotsb = x_n = -1$.

Rather surprisingly, this problem somewhat resembles intersecting set systems. 
If we want $\mathscr{F} \subset \mathscr{P}([n])$ to satisfy $A \cap B \neq \emptyset$ for all $A ,B \in \mathscr{F}$, then $\mathscr{F}$ can have at most $2^{n-1}$ members, since $\mathscr{F}$ can contain at most one of $C$ and $[n]\setminus C$ for each $C \subset [n]$. The bound can be attained in many different ways, for example by the family $\mathscr{F} = \left\{A\subset [n]\colon 1 \in A \right\}$.

What happens if we restrict ourselves to sets of a fixed size? In other words, what is the minimal number of non-negative $k$-wise sums, given $\sum \mathcal{X} \geq 0$? This question is essentially still open. The above example of $x_1 = n-1 $ and $x_2 = x_3 = \dotsb = x_n = -1$ gives $\binom{n-1}{k-1}$ non-negative sums. We can also consider the `mirror image' of that example: $x_1 = x_2 = \dotsb x_{n-1} = 1$ and $x_n = -n+1$ gives $\binom{n-1}{k}$ non-negative sums, which is less than $\binom{n-1}{k-1}$ for $n<2k$ and equals $\binom{n-1}{k-1}$ for $n=2k$.

This situation again parallels the behaviour of intersecting set systems. The Erd\H{o}s-Ko-Rado theorem \cite{EKR} states that for $n \geq 2k$ a family of sets $\mathscr{F}\subset [n]^{(k)}$ with the property that $A\cap B \neq \emptyset$ for all $A,B \in \mathscr{F}$ has at most $\binom{n-1}{k-1}$ elements. For $n>2k$ the bound is attained uniquely up to an isomorphism by $\mathscr{G} = \left\{A\in [n]^{(k)}\colon 1 \in A \right\}$.

The relation between non-negative $k$-sums and intersecting $k$-uniform set systems is rather subtle. To the best of our knowledge, no obvious way of translating one problem into another has been found so far. It is conceivable that there exists a min-max correspondence in the spirit of linear programming, for in the former problem we try to \emph{minimise} the number of non-negative $k$-sums, whereas in the latter we want to \emph{maximise} the size of an intersecting $k$-set system.

Let $A(n,k)$ to be the minimal number of non-negative $k$-sums over all possible choices of $\mathcal{X}\in \mathbbm{R}^{(n)}$ with $\sum \mathcal{X} \geq 0$. For what values of $n$ and $k$ do we have $A(n,k)=\binom{n-1}{k-1}$? 

This question was first asked by Bier and Manickam (see \cite{Bi},\cite{BM}), who investigated the so-called first distribution invariant of the Johnson-scheme in the 1980's. The following conjecture was made by Manickam and Mikl\'os~\cite{MM} in 1987, and, in a slightly different context, by Manickam and Singhi~\cite{MS} in 1988. 

\begin{conjecture}\label{conj: MS}
For all $n \geq 4k$ we have $A(n,k)=\binom{n-1}{k-1}$.
\end{conjecture}

Manickam and Mikl\'os \cite{MM} proved a number of results supporting this conjecture, including the following four assertions.

\begin{itemize}
\item $A(n,k)=\binom{n-1}{k-1}$ if $n$ is a multiple of $k$.
\item Conjecture~\ref{conj: MS} holds for $k=2$ and $k=3$.
\item $A(n,k)<\binom{n-1}{k-1}$ for some small values of $n$ like $n=3k+1$.
\item There is a function $f\colon \mathbb{N}\rightarrow \mathbb{N}$ such that for every $n\ge f(k)$ we have $A(n,k)=\binom{n-1}{k-1}$.

\end{itemize}

The first assertion is a direct corollary of the deep and powerful Baranyai partition theorem, proved in 1975 \cite{Bry}, which states that if $k|n$ then $[n]^{(k)}$ can be partitioned into blocks, each of which consists of $n/k$ pairwise disjoint $k$-sets. Since $\sum \mathcal{X}\geq 0$, each block must contain at least one non-negative $k$-set, whence we conclude that $A(n,k)\geq \frac{k}{n}\binom{n}{k} = \binom{n-1}{k-1}$. There are also proofs that avoid the usage Baranyai's theorem, see \cite{MM} for more details.

It is not hard to check by hand that $f(2)=6$. The case $k=3$, which needs a little more work, was settled by Manickam~\cite{Mani} and, independently and much later, by Marino and Chiaselotti~\cite{Marino}.

The `counterexample' for $n=3k+1$ is as follows: let $\mathcal{X}$ comprise number $3$, taken $3k-2$ times and number $-3k+2$, taken $3$ times. Since $3(k-1)+(-3k+2)=-1 < 0$, the non-negative $k$-wise sums are just those of the $3$'s, and there are $\binom{n-3}{k}$ of them. The inequality \[\binom{n-3}{k}<\binom{n-1}{k-1}\] can be re-written after some standard manipulations of binomial coefficients as \[(n-k)(n-k-1)(n-k-2)< k(n-1)(n-2).\] Substituting $n=3k+1$ we obtain \[(2k+1)\cdot 2k\cdot (2k-1)<k\cdot 3k(3k-1),\] which eventually simplifies to $(k-1)(k-2)>0$, proving that for $k>2$ this family has indeed fewer than $\binom{n-1}{k-1}$ non-negative $k$-sums. The above calculation was essentially based on the fact that $3^2>2^3$. Since this is no longer true for $3$ and $4$ in place of $2$ and $3$, the construction does not extend analogously to $4k+1$ and larger values of $n$.

Let $f(k)$ be the minimal $n_0$ such that $A(n,k)=\binom{n-1}{k-1}$ for all $n>n_0$. Conjecture~\ref{conj: MS} states that $f(k)\leq 4k$. The fourth of the listed assertions states that $f(k)$ is well-defined for every $k$. The best known bound so far was $f(k) \leq (k-1)(k^k+k^2)+k$, proved by Manickam and Mikl\'os in 1987 \cite{MM}; it has remained essentially unbeaten for the past 20 years. In 2003 Bhattacharya~\cite{Bh} gave a new proof of the existence of $f$, but his bound did not improve on the above. 

Another theorem proved by Manickam and Mikl\'os in~\cite{MM} states that $A(n,k)=\binom{n-1}{k-1}$ implies $A(n+k,k)=\binom{n+k-1}{k-1}$ and $A(cn,k)=\binom{cn-1}{k-1}$ for all $c$. As a consequence, if we can show that $A(n,k)=\binom{n-1}{k-1}$ for some $n$ coprime with $k$, we can deduce $f(k)\leq (k-1)n$. We think this fact might be very useful for proving Conjecture~\ref{conj: MS}, but we are not going to apply it here.

\section{Main result}

Our aim in this article is to establish a new bound on $f(k)$. We shall prove that $f(k)\leq k(4e\log k)^k = \exp(k\log \log k +O(k))$, which is a substantial improvement of the bound in \cite{MM}. But before doing this, we shall give a quick proof of $f(k) \leq 3k^{k+1}+k^3$. This is only slightly worse than the bound in~\cite{MM}, but it is derived using a completely new method. The ideas introduced here will be helpful later in the proof of our main result.

\begin{theorem}\label{thm:mansin}
If $n \geq 3k^{k+1}+k^3$, then the number of non-negative $k$-wise sums of real numbers, whose total sum is non-negative, is at least $\binom{n-1}{k-1}$.
\end{theorem}

\begin{proof}

From here on we assume without loss of generality that $x_1\geq x_2 \geq \dots \geq x_n$. If $x_1$ forms a non-negative sum with $k-1$ smallest elements of $\mathcal{X}$, then we obtain $\binom{n-1}{k-1}$ non-negative $k$-sums involving $x_1$, and we are done. So let us assume that 

\begin{equation}\label{order1} 
x_1 + \sum_{i=n-k+2}^{n}x_i < 0. 
\end{equation}

We claim the existence of at least $\binom{n-2k}{k-1}$ non-negative $k$-sums using neither $x_1$ nor one of $x_{n-k+2}, x_{n-k+3}, \dots x_n$. Since $\mathcal{X}\setminus\left\{x_1, x_{n-k+2}, x_{n-k+3}, \dotsc x_n \right\}$ has a non-negative total sum, we can throw away a few more negative members of $\mathcal{X}$ to obtain a collection, whose total sum is non-negative and whose size $m$ is a multiple of $k$ larger than $n-2k$. 

Since Conjecture~\ref{conj: MS} holds when $n$ is a multiple of $k$, the number of non-negative $k$-sums in such a collection is at least $\binom{m-1}{k-1}\geq\binom{n-2k}{k-1}$. Here we are assuming that at least $2k$ members of $\mathcal{X}$ are negative. This is a valid assumption, since otherwise there would be at least $\binom{n-2k}{k}$ non-negative $k$-sums, which is more than $\binom{n-1}{k-1}$ for $n \geq 3k^{k+1}+k^3$ (or in fact for $n>3k^2$, as can be shown by some more careful estimates).

Next we claim that the number of non-negative $k$-sums involving $x_1$ is at least $\binom{\left\lfloor n/k\right\rfloor}{k-1}$. More precisely, we claim that taking $\mathcal{Z}=\left\{x_2,x_3,\dots, x_{\left\lceil n/k\right\rceil}\right\}$ to be the collection of $\left\lfloor n/k \right\rfloor$ largest numbers in $\mathcal{X}\setminus \{x_1\}$ , any $k-1$ members of $\mathcal{Z}$ have a sum of at least $-x_1$ and therefore yield a non-negative sum when taken together with $x_1$. 

This follows from the fact that for every $j$ we have \[jx_1+(n-j)x_{j+1}\geq\sum\mathcal{X}\geq 0,\] thus 

\begin{equation}\label{order2}
x_1 \geq -x_{j+1}(\frac{n}{j}-1)=-x_{j+1}\frac{n-j}{j}.
\end{equation}

For $j=\lfloor n/k \rfloor$ this translates to \[x_1+\sum_{i=1}^{k-1}x_{a_i} \geq x_1+x_{n/k}(k-1)\geq 0,\] where the $x_{a_i}$ are any $k-1$ members of $\mathcal{Z}$.

Therefore, there must be at least $\binom{n-2k}{k-1}+\binom{\left\lfloor n/k\right\rfloor}{k-1}$ non-negative $k$-sums. The following calculation shows that for $n \geq 3k^{k+1}+k^3$ this is more than $\binom{n-1}{k-1}$, proving the theorem.

After applying the obvious estimates for binomial coefficients and multiplying through by $(k-1)!$ it remains to show that 
\[(n-3k)^{k-1}+(n/k-k)^{k-1}\geq n^{k-1}.
\]
From elementary probability we know that for $p, q > 0$ and $m\in \mathbbm{N}$ the values of monomials in the expansion of $(p+q)^m$ form a unimodal sequence, i.e.\ increase until the maximum is reached and decrease afterwards. In particular, if the first term in the expansion is larger than the second, the sequence of terms is monotone decreasing and so $(p-q)^m$ can be bounded below by the difference between the first and the second term. Let us check, if this condition is satisfied for the left hand side of the above inequality. We have

\[n^{k-1} > (k-1)n^{k-2}\cdot 3k,
\]
since $n > 3k(k-1)$ and

\[(n/k)^{k-1} > (k-1)(n/k)^{k-2}\cdot k, \] 
since $n > k^2(k-1)$. Therefore
\[(n-3k)^{k-1}+(n/k-k)^{k-1} > n^{k-1} - (k-1)n^{k-2}\cdot 3k + (n/k)^{k-1} - (k-1)(n/k)^{k-2}\cdot k 
\]
\[> n^{k-1} - n^{k-2}\cdot 3k^2 + \frac{n^{k-1}}{k^{k-1}} - \frac{n^{k-2}}{k^{k-4}}
\]
So it suffices to show that 
\[n^{k-1} - n^{k-2}\cdot 3k^2 + \frac{n^{k-1}}{k^{k-1}} - \frac{n^{k-2}}{k^{k-4}} \geq n^{k-1},
\]
which is equivalent to $n \geq 3k^{k+1}+k^3$. 
\end{proof}

Doing the estimates more carefully, the bound of $3k^{k+1}+k^3$ can be reduced by a constant factor. However, using a slightly different method we can do a lot better than this. Our main result is the following. 

\begin{theorem}\label{thm:mansin2}
For all $k>k(4e\log k)^k$ we have $A(n,k)=\binom{n-1}{k-1}$. In other words, $f(k)\leq k(4e\log k)^k$.
\end{theorem}

\begin{proof}
 
As in the proof of Theorem \ref{thm:mansin}, we may assume \eqref{order2}. Let us call an element $y$ in the set $\mathcal{Y}$ of real numbers \emph{central} if all $k$-sums in $\mathcal{Y}$ involving $y$ are non-negative. 

As before, if $x_1$ is central in $\mathcal{X}$, we are done. So let us assume that \eqref{order1} holds. 
Define $\mathcal{X}_2=\mathcal{X}\setminus \left \{x_1,x_{n-k+2},x_{n-k+3},\dots, x_n \right\}$; we know from \eqref{order1} that $\sum \mathcal{X}_2\geq 0$. Hence, by the same argument as in \eqref{order2}, we conclude that $x_2\geq - \frac{n-k-j}{j}x_{j+2}$ for all $j$.

Suppose that $x_2$ is central in $\mathcal{X}_2$. Then there are $\binom{n-k-1}{k-1}$ non-negative $k$-sums in $\mathcal{X}_2$ involving $x_2$. If we replace in each of them $x_2$ with $x_1$, we obtain just as many sums involving $x_1$ but not $x_2$. Therefore in total we obtain at least $2\binom{n-k-1}{k-1}$ non-negative sums. Using the first term -- second term estimate as in the proof of Theorem \ref{thm:mansin}, one can see that this is more than $\binom{n-1}{k-1}$ for $n>2k^2$. 

So let us assume that $x_2$ is not central in $\mathcal{X}_2$. Remove $x_2$ and the smallest $k-1$ members of $\mathcal{X}_2$ and call the resulting set $\mathcal{X}_3$; we know that $\sum \mathcal{X}_3>0$. Again, if $x_3$ is central in it, we obtain $3\binom{n-2k-1}{k-1}$, which is bigger than $\binom{n-1}{k-1}$ in the considered range, so we can remove $x_3$ and the smallest $k-1$ numbers. Repeat the procedure. 

We want to iterate the above argument $n/2k$ times. For what values $n$ is this possible? We need to make sure that $(p+1)\binom{n-kp-1}{k-1}>\binom{n-1}{k-1}$ for each $p$ between $1$ and $n/2k$. It suffices to prove the stronger statement that $(p+1)(n-k(p+1))^{k-1}>n^{k-1}$. Let us do it separately for $p<n/k^2$ and for $n/k^2 \leq p <n/2k$.

In the case $p<n/k^2$ we would like to estimate the left hand side using the first term -- second term method as in the proof of Theorem \ref{thm:mansin}. The condition $n^{k-1}>(k-1)k(p+1)n^{k-2}$ can be easily confirmed for $p<n/k^2$ and $n>k^3$. The estimate tells us that \[(p+1)(n-k(p+1))^{k-1}>(p+1)n^{k-1}-(p+1)^2k(k-1)n^{k-2}.\] By a straightforward manipulation of terms, the right hand side is greater than $n^{k-1}$ if and only if \[\frac{n}{k(k-1)}> p + 2 + \frac{1}{p}.\] The right hand side of this is at most $n/k^2 + 3$, which is less than the left hand side for $n>3k^2(k-1)$.

In the case $n/k^2 \leq p <n/2k$, we can estimate \[(p+1)(n-k(p+1))^{k-1}>\frac{n}{k^2}(n-k\frac{n}{2k})^{k-1}=\frac{n^k}{2^{k-1}k^2},\] which is greater than $n^{k-1}$ for $n>2^{k-1}k^2$.

So, suppose that for each $i$ between $1$ and $n/2k$, $x_i$ is not central in $\mathcal{X}_i$. Similarly to \eqref{order2}, for $x_{n/2k}$, the largest element of $\mathcal{X}_{n/2k}$, we obtain for each $1\leq j \leq n/2$ that \[x_{n/2k} \ge -\frac{n/2-j}{j}x_{j+n/2k}.\] Choosing $j=\frac{n}{2\log k}$, this translates to $x_{n/2k}>-(\log k - 1)x_{n/2 \log k + n/2k}$. Hence the sum of any $x_i$ with $1 \le i \le n/2k$ and $\log k -1 $ numbers $x_j$, each $j$ lying between $n/2k$ and $n/2k+n/2 \log k$, is positive. Thus we may pick $k/\log k$ numbers $x_i$ in the former range and $k-k/\log k$ numbers $x_j$ from the latter to obtain a non-negative $k$-sum. The total number of such sums will be \[\binom{n/2k}{k/\log k}\binom{n/2\log k}{k-k/\log k}.\] It remains to check that this is greater than $\binom{n-1}{k-1}$ for $n>k(4e\log k)^k$. Straightforward estimates yield

\[
\binom{n/2k}{k/\log k}
\geq \frac{n^{k/\log k}}{(4k)^{k/\log k}(k/\log k)!} \]

and

\[ \binom{n/2\log k}{k-k/\log k} \geq 
\frac{n^{k-k/\log k}}{(4 \log k)^{k-k/\log k}(k-k/\log k)!}.\]

Therefore, we obtain

\[ \binom{n/2k}{k/\log k}\binom{n/2\log k}{k-k/\log k} \geq \frac{n^k}{k!}\left(\frac{1}{4k}\right)^{k/\log k}\left(\frac{1}{4 \log k}\right)^{k-k/\log k} 
\]
\[ > \frac{n^k}{k!\cdot 4^k \cdot k^{k/ \log k}\cdot (\log k)^k} = \frac{n^k}{k!(4e\log k)^k}.
\]


The last expression is greater than $\frac{n^{k-1}}{(k-1)!}$, which is always greater than $\binom{n-1}{k-1}$, precisely for $n>k(4e\log k)^k =\exp(k \log \log k +O(k))$.

\end{proof}

\section*{Acknowledgements}
\label{sec:acknowledgements}

I would like to thank Dezs\H{o} Mikl\'os for drawing this problem to my attention.

\end{document}